\documentclass[14pt]{article}
\usepackage{amsmath,amssymb,amsthm,amsfonts}

\begin{document}
\bibliographystyle{plain}


\newtheorem{theor}{Theorem}[section] 
\newtheorem{prop}[theor]{Proposition} 
\newtheorem{cor}[theor]{Corollary}
\newtheorem{lemma}[theor]{Lemma}
\newtheorem{sublem}[theor]{Sublemma}
\newtheorem{defin}[theor]{Definition}
\newtheorem{conj}[theor]{Conjecture}


\newcommand{\beginProof}{\par{\bf Proof: }}
\newcommand{\ProofEnd}{{\bf Q.E.D.}}

\newcommand{\kkk}[1]{{\Large\bf#1}}
\newcommand{\nlabel}[1]{\label{#1}}


\newcommand{\mF}{{\mathbb F}}
\newcommand{\mG}{{\mathbb G}}
\newcommand{\mP}{{\mathbb P}}
\newcommand{\mQ}{{\mathbb Q}}
\newcommand{ \mR}{{\mathbb R}}
\newcommand{ \mZ}{{\mathbb Z}}
\newcommand{\mC}{{\mathbb C}}
\newcommand{\mN}{{{\mathbb Z}_{\geqslant 1}}}
\newcommand{\Qb}{\overline{\mathbb Q}}


\newcommand{\Div}{{\rm Div}}
\newcommand{\mtr}[1]{\overline{#1}}
\newcommand{\Stab}{{\rm Stab}}
\newcommand{\Irr}{{\rm Irr}}
\newcommand{\Zar}{{\rm Zar}}
\newcommand{\Sup}{{\rm Sup}}
\newcommand{\Tor}{{\rm Tor}}
\newcommand{\Gal}{{\rm Gal}}
\newcommand{\Spec}{{\rm Spec}}


\newcommand{\ra}{\rightarrow}
\newcommand{\without}{\setminus}


\author{Damian R\"ossler
\footnote{CNRS, Centre de Math\'ematiques de Jussieu,  
     Universit\'e Paris 7 Denis Diderot,  
     Case Postale 7012,  
     2, place Jussieu,  
     F-75251 Paris Cedex 05, France,   
     E-mail : dcr@math.jussieu.fr}}
\title{An afterthought on the generalized Mordell-Lang conjecture}
\maketitle
\begin{abstract}
The generalized Mordell-Lang conjecture (GML) is the statement that the irreducible components of 
the Zariski closure 
of a subset of a group of finite rank inside a semi-abelian variety are translates 
of closed algebraic subgroups. 
In \cite{McQuillan}, M. McQuillan gave a proof of this statement. We revisit his proof, indicating 
some simplifications. This text contains a complete elementary proof 
of the fact that (GML) for groups of torsion points (=\ generalized Manin-Mumford conjecture), together with (GML) 
for finitely generated groups  imply the full 
generalized Mordell-Lang conjecture.
\end{abstract}
\date
\begin{center}
\end{center}


\parindent=0pt
\parskip=5pt

\section{Introduction}

Let $A$ be a semi-abelian variety over $\Qb$ (cf. beginning of 
Section 2 for the definition of a semi-abelian variety). We shall call a closed 
reduced subscheme of $A$ {\it linear } if its irreducible components 
are translates of closed subgroup schemes of $A$ by 
points of $A(\Qb)$.  
Let $\Gamma$ be a finitely generated subgroup of 
$A(\Qb)$ and define  
$\Div(\Gamma):=\{a\in A(\Qb) | \exists n\in\mZ_{\geqslant 1}: n\cdot a\in\Gamma\}$. 
Let $X$ be a closed reduced subscheme of $A$. 
Consider 
the following statement : 

{\it the variety $\Zar(X\cap\Div(\Gamma))$ is linear} (*).

The generalized Mordell-Lang  conjecture (GML) is the statement that 
(*) holds for any data $A,\Gamma, X$ as above. 
The statement (GML) with the supplementary requirement
 that $\Gamma=0$ shall be referred to as (MM). 
The statement (GML) with $\Div(\Gamma)$ replaced by $\Gamma$ in (*) shall be referred to 
as (ML). 

The statement (ML) was first proven by Vojta (who built on Faltings work) in \cite{Vojta}. 
The statement (MM) was first proven by Hindry (who built on work of Serre and Ribet) in \cite{Hindry}.  Finally, McQuillan (who built 
on the work of the previous) proved (GML) in \cite{McQuillan}. 

The structure of McQuillan's proof of (GML) has three key inputs:
(1) the statement (ML), (2) an extension of (MM) to families of varieties 
 and (3) the 
Kummer theory of abelian varieties.

In this text, we shall indicate some simplifications of this proof. 
More precisely, we show the following. First, that once  (MM) is granted, 
a variation of (2) sufficient for the purposes of the proof is contained in an automatic uniformity principle proved by Hrushovski. 
See Lemma \ref{lemmaunif} 
for the statement of this automatic uniformity principle and a reference for the proof, which uses nothing 
more than the compactness theorem of first order logic. 
Secondly, 
we show that one can replace the Kummer theory 
of abelian varieties (3) by an elementary geometrical argument.
The core of the simplified proof is thus an elementary proof of the following 
statement:

{\it if {\rm (ML)} and {\rm (MM)} hold then {\rm (GML)} holds} (**)

and the proof of (GML) is then obtained by combining (**) with 
the statements (MM) and (ML), which are known to be true by the work 
of Hindry and Vojta. 

We stress that our proof of (**) is independent of 
the truth or techniques of proof of either (ML) or (MM). 

Notice that yet another proof of (GML) was given 
by Hrushovski in \cite[Par. 6.5]{Hrushovski}. His proof builds on 
a generalisation of his model-theoretic proof of (MM) (which 
is based on the dichotomy theorem of the theory of 
difference fields)  and on (ML). It also avoids the Kummer theory of abelian varieties 
but it apparently doesn't lead to a proof of (**).
Finally, we want to remark that a deep Galois-theoretic 
result of Serre (which makes an earlier statement of Bogomolov 
uniform), which is used in McQuillan's proof of (GML) as well 
as in Hindry's proof of (MM) (see \cite[Lemme 12]{Hindry}), was never published. 
Now different proofs of (MM), which do not rely on Serre's result, 
were given by Hrushovski in \cite{Hrushovski} and by Pink-R\"ossler in \cite{PR2}. 
Our proof of (**) thus leads to a proof of (GML) which is independent of Serre's result.

The structure of the article is as follows. Section 3 contains 
the proof of (**) and section 2 recalls 
the various facts from the theory of semi-abelian varieties that we shall 
need in Section 3. The reader is encouraged to 
start with Section 3 and refer to Section 2 as necessary.

{\bf Basic notational conventions.}  
A (closed) subvariety of a scheme $S$ is a (closed) reduced subscheme of $S$.
If $X$ is a closed subvariety of a $\Qb$-group scheme $A$, we shall write $\Stab(X)$ for 
the stabilizer of $X$ in $A$, which is a 
closed group subscheme of $A$ such that $\Stab(X)(\Qb)=\{a\in A(\Qb) | X+a=X\}$. 
If $H$ is a commutative group, we write $\Tor(H)\subseteq H$ for the subgroup 
consisting of the elements of finite order in $H$. If $T$ is a noetherian topological space, denote 
by $\Irr(T)$ the set of its irreducible components. 

\section{Preliminaries}

A {\it semi-abelian variety} $A$ over an algebraically closed field $k$ is 
by definition a commutative 
group scheme over $k$ with the following properties: 
it has a closed subgroup scheme $G$ which is 
isomorphic to a product of finitely many multiplicative groups over $k$ and there 
exists an abelian variety $B$  over $k$ and a surjective morphism 
 $\pi:A\ra B$, which is a morphism of group schemes over $k$ and whose kernel is $G$. 

In the next lemma, let $A$ be a semi-abelian variety as in the previous definition.
\begin{lemma}
Let $n\in\mZ_{\geqslant 1}$. 
The multiplication by $n$ morphism $[n]_A:A\ra A$ is quasi-finite. 
\label{lemmaqf}
\end{lemma}
\beginProof
we must prove that the fibers of 
$[n]_A$ have finitely many points or equivalently that 
they are of dimension $0$. Moreover, since 
the function $\dim([n]_{A,a})$ (=\ dimension of the fiber of $[n]_A$ over $a$) is 
a constructible function of $a\in A$ (see \cite[Ex. 3.22, chap. II]{Hartshorne}), 
it is sufficient to prove that the fibers of $[n]_A$ over 
closed points are finite. A fiber of $[n]_A$ 
over a closed point can be identified with the fiber $\ker\ [n]_A$ of $[n]_A$ 
over $0\in A$. The scheme $\ker\ [n]_A$ is naturally fibered over 
$\ker\ [n]_B$. The scheme $\ker\ [n]_B$ consists of finitely many 
closed points because multiplication by $n$ in 
$B$ is a finite morphism, as $B$ is an abelian variety (see \cite[Prop. 8.1 (d)]{Milne}. It will thus be sufficient to prove 
that the fibers of the morphism 
$\ker\ [n]_A\ra\ker\ [n]_B$ are finite and furthermore each of these fibers can identified with the fiber of 
$\ker\ [n]_A\ra\ker\ [n]_B$ over $0$. By construction this fiber is the closed subscheme 
 $\ker\ [n]_A\times_A G=\ker\ [n]_G$ of $A$. To prove that 
$\ker\ [n]_G$ has finitely many closed points, choose 
an identification $G\simeq G_m^\rho$ of $G$ with a product of $\rho$ 
multiplicative groups over $k$. The closed points 
of $\ker\ [n]_G$ then correspond to $\rho$-tuples of 
$n$-th roots of unity in $k$. This set is finite and this concludes the proof of the lemma.
\ProofEnd

Let now $A$ be a semi-abelian variety over $\Qb$. 

\begin{theor}[Kawamata-Abramovich]
Let $X$ be a closed 
subvariety of $A$. 
The union $Z(X)$ of the irreducible linear subvarieties of positive dimension of $X$ is Zariski closed.
The stabilizer $\Stab(X)$ of $X$ is finite if and only if the complement of $Z(X)$ in $X$ is not empty. 
\label{theorkaabr}
\end{theor}
For the proof see \cite[Th. 1 \& 2]{Abramovich} .

Let $Y$ be a variety over $\Qb$ and let $W\hookrightarrow A\times_{\Qb} Y$ be a closed subvariety. 

\begin{lemma}[Hrushovski]
If {\rm (MM)} holds then the quantity
$$
\Sup\{\#\Irr\big(\Zar(W_y\cap\Tor(A(\Qb)))\big)\}_ {y\in Y(\Qb)}
$$
is finite.
\label{lemmaunif}
\end{lemma}
Notice that (MM) predicts that the irreducible components of 
$\Zar(W_y\cap\Tor(A(\Qb)))$ are linear for each $W_y$, $y\in Y(\Qb)$. 
In words, the content of the lemma is that if this is the case, then the number of these irreducible 
components can be bounded independently of $y\in Y(\Qb)$. 
A self-contained proof of Lemma \ref{lemmaunif} can be found in 
\cite[Postscript, Lemma 1.3.2, p. 52-53]{Hrushovski}. For an extension of 
Lemma \ref{lemmaunif}, see \cite{Scanlon}. 

Suppose now that $A$ has a model $A_0$ over a number field 
$K$. 

In \cite[Par. 1.2 and Prop. 2]{Serre} it is shown that 
there exists a variety $\mtr{A}_0$ projective over $K$ and 
an open immersion $A_0\hookrightarrow\mtr{A}_0$ such that 
for all $n\in\mZ_{\geqslant 1}$ 
the multiplication by $n$ morphism $[n]_{A_0}:A_0\ra A_0$ extends to a $K$-morphism 
$[n]_{\mtr{A}_0}:\mtr{A}_0\ra \mtr{A}_0$. Furthermore, it is shown in \cite[Prop. 3]{Serre} that 
the corresponding diagram
\begin{equation}
 \begin{matrix}
    A_0 & \hookrightarrow& \mtr{A}_0\\
     \downarrow [n]_{A_0}& & \downarrow [n]_{\mtr{A}_0}\\
    A_0 & \hookrightarrow& \mtr{A}_0
    \end{matrix}
\nonumber
\end{equation}
is then cartesian. 

Let $\Gamma$ be a finitely generated subgroup of 
$A_0(K)$. 
\begin{lemma}[McQuillan]
The group generated by the set $\Div(\Gamma)\cap A_0(K)$ is 
finitely generated.
\label{lemmamcq}
\end{lemma}
Lemma \ref{lemmamcq} is McQuillan's Lemma 3.1.3 in \cite{McQuillan}. 
As the proof given there is somewhat sketchy, we shall provide 
a proof of Lemma \ref{lemmamcq}. 
\beginProof
let $B$ be an abelian variety over $\Qb$ and $\pi:A\ra B$ be a $\Qb$-morphism 
whose kernel $G$ is isomorphic to a product of tori over $\Qb$. These data 
exist since $A$ is a semi-abelian variety. 
Notice that for the purposes of the proof we may enlarge the field of definition $K$ of $A$ if necessary, 
since that operation will also enlarge the set $\Div(\Gamma)\cap A_0(K)$. 
Hence we may assume that $B$ (resp. $G$) has a model $B_0$ (resp $G_0$)
over $K$ and that $\pi$ has a model $\pi_0$ over $K$. Furthermore, 
we may assume that the isomorphism of $G$ with a product of 
tori descends to a  $K$-isomorphism of $G_0$ with a product of split tori over $K$.
Now fix a compactification $\mtr{A}_0$ of $A_0$ over $K$ as above. 
We consider the following situation. 
The symbol  $V$ refers to an open subscheme of the spectrum $\Spec\ {\cal O}_K$ of the ring of integers 
of $K$ and 
 ${\cal A}_0$ is a semi-abelian scheme over $V$, which is a model of $A_0$. 
 The symbol 
$\mtr{\cal A}_0$ refers to a projective model over $V$ of $\mtr{A}_0$ and we suppose given 
an open immersion ${\cal A}_0\hookrightarrow\mtr{\cal A}_0$, which is a model of 
the open immersion $A_0\hookrightarrow\mtr{A}_0$. We also suppose that 
the multiplication by $n$ morphism $[n]_{{\cal A}_0}$ on ${\cal A}_0$ extends to a $V$-morphism 
$[n]_{\mtr{\cal A}_0}:\mtr{\cal A}_0\ra \mtr{\cal A}_0$ and we suppose that 
the corresponding diagram 
\begin{equation}
 \begin{matrix}
    {\cal A}_0 & \hookrightarrow& \mtr{\cal A}_0\\
     \downarrow [n]_{{\cal A}_0}& & \downarrow [n]_{\mtr{\cal A}_0}\\
    {\cal A}_0 & \hookrightarrow& \mtr{\cal A}_0
    \end{matrix}
\nonumber
\end{equation}
is cartesian. 
Let ${\cal B}_0$ (resp. ${\cal G}_0$) be a model of $B_0$ (resp. $G_0$)
over $V$ and let $\widetilde{\pi}_0$ be a model of $\pi_0$ over $V$. Furthermore, 
we assume that the $K$-isomorphism of $G_0$ with a product of split 
tori extends to a $V$-isomorphism of ${\cal G}_0$ with a product of split tori over $V$.
We leave it to the reader to show that there are objects $V$, ${\cal A}_0$ etc. 
satisfying the described conditions.

The morphism $[n]_{{\cal A}_0}$ is then 
proper, because $[n]_{\mtr{\cal A}_0}$ is proper (as $\mtr{\cal A}_0$ is proper 
over $V$)  and properness is invariant 
under base change. The morphism $[n]_{{\cal A}_0}$ is therefore 
finite, as it is quasi-finite by \ref{lemmaqf} (applied to each fiber of 
${\cal A}_0$ over $V$). 
We may suppose without restriction of generality that 
$\Gamma$ lies in the image of ${\cal A}_0(V)$ in $A_0(K)$; indeed this condition 
will always be fulfilled after possibly removing a finite number of closed points from 
$V$. 
Let $a\in\Div(\Gamma)\cap A_0(K)$. Choose an $n\in\mZ_{\geqslant 1}$ such that 
$n\cdot a\in\Gamma$. 
Let $E$ be the image of 
the section $V\ra{\cal A}_0$ 
corresponding to $n\cdot a$. 
Consider the reduced 
irreducible component $C$ of 
$[n]_{{\cal A}_0}^*E$ containing the image of $a$. 
The image of $a$ is the generic point of 
$C$ and by assumption the natural morphism $C\ra E$ identifies 
the function fields of $C$ and $E$.  Furthermore, 
the morphism $C\ra E$ is finite. 
Now let $R_0$ be the ring underlying the affine scheme $V$. 
In view of the above, we can write $C=\Spec\ R$, where 
$R$ is a domain and the morphism $C\ra E$ identifies 
$R$ with an integral extension of $R_0$ inside the integral closure of  
$R_0$ in its own field of fractions. As $R_0$ is integrally closed (it is even a Dedekind ring) 
$C\ra E$ is an isomorphism. Hence 
$a\in{\cal A}_0(V)$. Thus we only have to show that 
${\cal A}_0(V)$ is finitely generated. 
This follows from the fact that 
${\cal B}_0(V)$ is finitely generated by 
the Mordell-Weil theorem applied to $B_0$ and the fact that 
${\cal G}_0(V)$ is finitely generated by  
the generalized Dirichlet unit theorem (see \cite[chap. V, par. 1]{Lang}).
\ProofEnd

\begin{lemma}
Let $C>0$. 
The set $\{a\in \Tor(A(\Qb)) | [K(a):K]<C\}$ is finite.
\label{lemmator}
\end{lemma}
\beginProof
if $A_0$ is an abelian variety over $K$ then this follows 
from the fact that the N\'eron-Tate height of torsion points vanishes 
and from Northcott's theorem. We leave the general case 
as an exercise for the reader.
\ProofEnd

\section{Proof of (**)}

In this section we shall prove (GML) using the results listed 
in Section 2 as well as (ML) and (MM). 

So we set off to prove (*). We may assume without loss of generality that 
$X$ is irreducible. We may also suppose that $\Stab(X)=0$. 

To see the latter, consider the closed subvariety $X/\Stab(X)$ of the 
quotient variety $A/\Stab(X)$. 
The image of the group $\Div(\Gamma)$ in $(A/\Stab(X))(\Qb)$ lies 
inside the group $\Div(\Gamma_1)$, where $\Gamma_1$ is the image of 
$\Gamma$ and the image 
of $\Div(\Gamma)$ is dense in $X/\Stab(X)$. So the assumptions of $(*)$ hold 
for $A/\Stab(X)$, $\Gamma_1$ and $X/\Stab(X)$. Furthermore, by construction 
$\Stab(X/\Stab(X))=0$. Now if (*) holds in this situation, $X/\Stab(X)$ is the translate 
of a connected closed group subscheme of $A/\Stab(X)$. By Theorem 
\ref{theorkaabr}, $X/\Stab(X)$ must therefore be a closed point. This 
implies that $X$ is a translate of $\Stab(X)$, thus proving (*) 
for $A$, $\Gamma$ and $X$. It is thus sufficient to prove 
(*) for $A/\Stab(X)$, $\Gamma_1$ and $X/\Stab(X)$ and we may thus 
replace $A$ by $A/\Stab(X)$, $\Gamma$ by $\Gamma_1$, 
$X$ by $X/\Stab(X)$. We then have $\Stab(X)=0$.

We may also assume without loss of generality that 
$A$ (resp. $X$) has a model $A_0$ (resp. $X_0$) over 
a number field $K$ such that $\Gamma\subseteq A_0(K)$ and such 
that the immersion $X\ra A$ has a model over $K$ as an immersion 
$X_0\ra A_0$. 

Let $U$ be the complement in $X$ of the union of the irreducible linear subvarieties of positive dimension of $X$. With a view to obtain a contradiction, we shall suppose that $U\not=\emptyset$. 
This implies that $U$ is dense in $X$ and in particular that $U\cap\Div(\Gamma)$ is dense in $X$. 
Let $a\in U\cap\Div(\Gamma)$ and let $\sigma\in\Gal(\Qb | K)$. By maximality 
$\sigma(U)\subseteq U$ and we thus have $\sigma(a)-a\in U-a$.
Now by definition, there exists $n=n(a)\in\mN$ such that $n\cdot a\in\Gamma\subseteq
A_0(K)$. We calculate
$$
n\cdot(\sigma(a)-a)=\sigma(n\cdot a)-n\cdot a=n\cdot a-n\cdot a=0.
$$
Thus $\sigma(a)-a\in\Tor(A(\Qb))\cap U-a$. The statement (MM) implies that 
$\Tor(A(\Qb))\cap U-a$ is finite and using Theorem \ref{lemmaunif}, we see 
that $\#(\Tor(A(\Qb))\cap U-a)<C$ for some constant $C\in\mN$, which is independent of $a$. 
Now this implies that $\#\{\tau(a)|\tau\in\Gal(\Qb | K )\}=\#\{\tau(a)-a|\tau\in\Gal(\Qb | K )\}<C$. 

A consequence of this conclusion is reached by McQuillan at the beginning of 
Par. 3.3 (p. 157) of \cite{McQuillan} using Theorem 3.2.2 of that article.

By Galois theory, we thus have 
$[K(a):K]<C$. We deduce from this last inequality that 
$$
[K(\sigma(a)-a):K]\leqslant [K(a,\sigma(a)):K]<C^2.
$$
By Lemma \ref{lemmator}, we see that this implies that 
$\sigma(a)-a\in T$, where $T\subseteq A(\mQ)$ is a finite 
set, which is dependent on $C$ but independent of either $a$ or $\sigma$. For each $b\in 
A(\mQ)\without A_0(K)$, choose 
$\sigma_b\in\Gal(\Qb | K)$ such that 
$\sigma_b(b)\not=b$. Suppose that 
the set $\{b\in U\cap\Div(\Gamma) | b\in A(\mQ)\without A_0(K)\}$ 
is dense in $X$. By the above, if this holds, there exists $t_0\in T$, $t_0\not=0$ such that the set 
$$
\{b\in U\cap\Div(\Gamma) | b\in A(\mQ)\without A_0(K),\ \sigma_b(b)-b=t_0\}
$$
is dense in $X$. Since $\sigma_b(b)\in U$ for all 
$b\in U$ such that $b\in A(\mQ)\without A_0(K)$, we see that this implies that 
$t_0\in\Stab(X)(\Qb)$. But $t_0\not=0$ so this contradicts 
our hypothesis that $\Stab(X)=0$. Thus we deduce that the set 
$\{b\in U\cap\Div(\Gamma) | b\in A(\mQ)\without A_0(K)\}$ is not dense in 
$X$ and thus the set $\{b\in U\cap\Div(\Gamma) | b\in A_0(K)\}$ 
must be dense in $X$. By Lemma \ref{lemmamcq}, 
the elements of this set are contained in a finitely generated group, 
so this contradicts (ML). So we have obtained a contradiction to our initial assumption 
that $U$ is not empty. Hence $U=\emptyset$ and Theorem \ref{theorkaabr} then shows that 
$X$ is a point. We have thus proven (GML) for $X$. 

The second part of our argument, based on Lemma \ref{lemmator}, is 
replaced by an argument involving the Kummer theories of abelian 
varieties and tori in Par 3.3 of McQuillan's article \cite{McQuillan}.

\end{document}